\documentclass[a4paper, 12pt]{amsart}

\newcommand{\version}{30/06/24 Anton}

\newcommand{\mytitle}[2][]{%
  \gdef\mylongtitle{#2}
  \gdef\myshorttitle{#1}
  \gdef\mypdftitle{#1}
}
\mytitle[Exponential no generic derivations]{Exponential Fields: Lack of Generic Derivations}
\newcommand*{\mykeywords}{Exponential fields, derivations, companionable theory}

\title[\myshorttitle{} v.\version]{\mylongtitle\\
}

\usepackage{lmodern}

\usepackage{microtype}

\author[A. Fornasiero]{Antongiulio Fornasiero}
\address{Universit\`a di Firenze}
    
\email{antongiulio.fornasiero@unifi.it} 
\urladdr{https://sites.google.com/site/antongiuliofornasiero/}

\author[G. Terzo]{Giuseppina Terzo}
\address{Universit\`a degli Studi di Napoli "Federico II"}
    
\email{giuseppina.terzo@unina.it} 

\date{04 June 2024}

\usepackage[svgnames]{xcolor}
\usepackage{ifpdf}
\ifpdf
\usepackage[pdftex, linktocpage=true, 
pdfborder={0 0 0}, plainpages=false,%
pdfpagelabels, pdfdisplaydoctitle=true,%
pdfstartview={FitV},
colorlinks=true,
citecolor=ForestGreen
]{hyperref}
\hypersetup{%
  pdfauthor={A. Fornasiero},
  pdftitle={\mypdftitle{} v.\version},
  pdfkeywords={\mykeywords},
  pdfsubject={First-order expansions of exponential fields do not have generic derivations}
}
\else
\usepackage[plainpages=false,linktocpage=true,pdfpagelabels, colorlinks=false
]{hyperref}
\fi

\usepackage[ocgcolorlinks]{ocgx2} 

\usepackage{amsmath, amsthm, amssymb}
\usepackage{mathtools}
\usepackage{xspace}

\usepackage{comment}

\usepackage{enumitem}
\setlist{
         left = 0pt}
\setlist[enumerate,1]{label = (\arabic*), 
         ref =   (\arabic*),
         font = \upshape}

\usepackage[T1]{fontenc}

\usepackage[english]{babel}
\usepackage[babel]{csquotes}
\usepackage[backend=biber, style=alphabetic, sorting=nyt,
 backref = true, backrefstyle = three,
 isbn = false, url = false
 ]{biblatex}
\AtEveryBibitem{%
  \clearlist{language}
  \clearfield{urlyear}
  \clearfield{urlmonth}
} 

\newcommand{\Td}{T^{\delta}}
\newcommand{\Tmorley}{T_{M}}
\newcommand{\Tmorleyd}{T_{M}^{\delta}}

\newcommand{\Ld}{L^{\delta}}
\newcommand{\Lmorley}{L_{M}}
\DeclareMathOperator{\ACFz}{ACF_0}
\DeclareMathOperator{\DCFz}{\mathrm{DCF}_0}
\newcommand{\ECM}{\mathcal E}

\newcommand*{\structure}[1]{\mathbb{#1}}
\newcommand{\K}{\structure{K}}
\newcommand{\Kdag}{\K^{\dagger}}
\newcommand{\Kd}{\K^\delta}
\DeclareMathOperator{\acld}{acl^{\delta}}

\newcommand*{\set}[1]{\{#1\}}

\newcommand*{\card}[1]{\lvert#1\rvert}
\newcommand{\N}{\mathbb{N}}
\newcommand{\Z}{\mathbb{Z}}
\newcommand{\R}{\mathbb{R}}
\newcommand{\C}{\mathbb{C}}
\newcommand{\Q}{\mathbb{Q}}

\newcommand*{\tuple}[1]{\langle #1 \rangle}

\newcommand{\bv}{\bar b}

\newcommand{\x}{\bar x}
\newcommand{\y}{\bar y}

\def\Ind#1#2{#1\setbox0=\hbox{$#1x$}\kern\wd0\hbox to
  0pt{\hss$#1\mid$\hss}\lower.9\ht0\hbox to 0pt{\hss$#1\smile$\hss}\kern\wd0}

\newcommand{\Cone}{\mathcal C^1}

\newcommand{\Wlog}{W.l.o.g\mbox{.}\xspace}

\newcommand{\eg}{e.g\mbox{.}\xspace}

\newtheorem{mainthm}{Theorem}

\newtheorem{lemma}{Lemma}[section]
\newtheorem{thm}[lemma]{Theorem}
\newtheorem{corollary}[lemma]{Corollary}

\newtheorem{open problem}[lemma]{Open problem}

\newtheorem*{fact*}{Fact}

\theoremstyle{remark}

\newtheorem{claim}{Claim}
\newtheorem*{claim*}{Claim}

\theoremstyle{definition}
\newtheorem{definition}[lemma]{Definition}

\newtheorem{remark}[lemma]{Remark}
\newtheorem{final remark}[lemma]{Final remark}

\begin{document}

\begin{abstract}
We investigate
the existence of ``generic derivations'' in exponential fields.
We show that exponential fields without additional compatibility conditions
between derivation and exponentiation cannot support a generic derivation.
\end{abstract}

\maketitle


\makeatletter
\renewcommand\@makefnmark%
   {\normalfont(\@textsuperscript{\normalfont\@thefnmark})}
\renewcommand\@makefntext[1]%
   {\noindent\makebox[1.8em][r]{\@makefnmark\ }#1}
\makeatother

\section{Introduction}
Let $L  \supseteq \{+, -, \cdot, 0, 1 \}$ be an expansion of the language of rings, 
and $T$ be an $L$-theory expanding the theory of fields of characteristic~$0$.
Let $\delta$ be a new unary function symbol, and $\Td$ be the $\Ld$ expansion of $T$
such that $\delta$ is a derivation; that is, $\delta$ satisfies the following axioms:
\[
\delta(x + y) = \delta x + \delta y, \quad
\delta(xy) = x \delta y + y \delta x.
\]
We denote by $\Tmorley$ the Morleyzation of $T$ (in the language $\Lmorley$).
We say that $T$ admits \emph{generic derivations} if $\Tmorleyd$ is companionable,
that is if it has a model
companion. 
The notion of model companion (see Definition~\ref{def:MC})
was introduced by Robinson;
for details on it and on Morleyzation, see \cite{Sacks}.
We will consider only the case when $T$ expands a field of characteristic~$0$.

Robinson \cite{Robinson} proved that $\ACFz$ admits generic derivations, and Blum gave an
explicit axiomatization for the theory $\DCFz$ of differentially closed fields
(see~\cite{Sacks}).

In recent years, mathematicians have shown that many other theories $T$ admit
generic derivations. 
The most general result to date is in~\cite{FT}, where we showed that
if $T$ is algebraically bounded, then it admits generic derivations; for a survey
of previous results, see the same article. 
There are also examples of theories
which do not admit a model companion, as the theory of groups, the theory of modules over certain rings
(for details see~\cite{Eklof}) and the theory of exponential fields 
(see~\cite{kirby}).

In this paper we focus on the case when $T$ expands a field equipped  with a
non-constant exponential function.
More precisely, we assume that $L$ includes an unary function symbol~$E$, 
and $T$ proves that
\[
E(x + y) = E(x) \cdot E(y), \qquad
\exists z\, E(z) \neq 1.
\]

We aim to prove that under these conditions $T$ does not admit generic derivations.

\smallskip

The existence of generic derivations has been established for some specific
theories~$T$,
but these results rely on the additional assumption that $\delta$ is an E-derivation.
In other words, the theory
$\Td$ also includes the extra axiom:
\[
\delta E(x) = E(x) \delta x.
\]
Furthermore, these results often require additional ``compatibility conditions''; for more details
see \cites{FK, Point, ADH}.

In this paper we are considering the case when no additional compatibility
conditions are imposed on the derivation~$\delta$: in particular, we don't require it
to be an E-derivation nor
continuous under some possible topology on models of~$T$.
The main result is the following:

\begin{mainthm}\label{thm:exp}
If $T$ expands a non-trivial exponential field, then $T$ does not admit generic derivations (see Theorem \ref{generale}).
\end{mainthm}

So we have as consequences that $(\mathbb C, \exp),$ $(\mathbb B, E(x))$ 
where $\mathbb B$ is a Zilber field (such fields are introduced by Zilber in
\cite{zilber}), and $(\R, \exp)$ do not admit generic derivation.

\smallskip

Using similar techniques we obtain analogues results for some o-minimal theories. In particular we prove that $(\R, \exp_{\upharpoonright[-1, 1]})$ and $(\R, \sin_{\upharpoonright [0, 2\pi]})$ does not admit generic derivations; for details see Sections~\ref{altricasi} and~\ref{caso2}.

\subsection*{Acknowledgements}
The authors thank Elliot Kaplan, Noa Lavi, Omar León Sanchez,  and Alex Usvyatsov
for useful discussions on the topic.

\section{Exponential function}
In this section we give the proof of the main result, i.e.\ Theorem~\ref{thm:exp}.  
We suppose that $T$ expands a non-trivial exponential field. \Wlog, we may assume that $T$ has quantifier elimination and we denote by $\K \models T$ and $\Kd  \models \Td$, both with domain~$K$.\\
The idea of the proof is to show that every existentially closed model $\Kd$ of
$\Td$ defines the set~$\Q$. 
Since a countable set cannot be definable in an  
$\aleph_{1}$-saturated model of a first-order theory, we obtain that the theory $\Td$ does not have model companion.

For completeness we recall some basic notions:

\begin{definition}
A model $M \models T$ is existentially closed if, for every quantifier free formula
$\phi(\x, \y)$ and all $\overline a$ in $M$, 
if there is an extension $M \subset B$ such that $B \models T$ and $B \models \exists \phi(\x, \overline
a)$ then 
\[M \models \exists \phi(\x, \overline a).\]
\end{definition}


\begin{definition}\label{def:MC}
 $\Td$ has a model companion if the class $\ECM$
of existentially closed
models of $\Td$ is elementary, and in that case the model companion 
is the
theory of $\ECM$.
\end{definition}

To lay the groundwork for the main result of this section, we first prove a key lemma. This lemma has broad applicability and will be instrumental in subsequent sections. Therefore, we present it in its general form.

\begin{lemma}\label{zariskidense}
Let $\K^{\delta} \in \ECM$.
Let $X \subset K^2$ be $\K$-definable, and consider $\delta_{\upharpoonright X} : X \rightarrow K^2$. 
We have that $\delta_{\upharpoonright X}$ is a surjective iff $X$ is Zariski
dense. (The only assumption on $T$ is that it is an expansion of fields of
characteristic~$0$, and of course $\delta$ is a  derivation). 
\end{lemma}

\begin{proof} 
$(\Rightarrow)$
For this direction, we do not need to assume that $\K^{\delta}  \in \ECM$
(but only that $\K^{\delta}\models \Td$).
It is equivalent to prove that, if $X$  is not Zariski dense, then the map $\delta_{\upharpoonright X}$ is not surjective.

Let ${\Kdag}^{\delta^{\dag}} \succ \K^{\delta}$ be $\card{K}^{+}$-saturated.
We remind that there is an induced pregeometry $\acld$ on $\Kdag$ (see \eg
\cite{FK}).
Let $a,b \in \Kdag$ be $\delta$-algebraically independent over~$K$.
Assume, by contradiction, that  $\delta$ is surjective, and therefore there exists $\tuple{u,v} \in X^{\dag}$ (the
interpretation of $X$ inside $\K^{\dag}$) such
that $\delta^{\dag} u = a$ and $\delta^{\dag} v = b$.
Thus, $u, v$ must be $\delta$-algebraically independent over $K$, but since
$\tuple{u,v} \in X^{\dag}$ they are
algebraically dependent over~$K$, we have a contradiction.

\medskip

$(\Leftarrow)$
Let $a,b \in K$; we want to prove that there exists $\tuple{u,v} \in X$ such
that $\delta u = a$ and $\delta v = b$.
Let $\Kdag \succ \K$ be a $\card{K}^{+}$-saturated $L$-structure (no derivation on
it yet).
Let $\tuple{u,v} \in X^{\dag}$ be algebraically independent over~$K$.
Let $\rho$ be a derivation on $\Kdag$ extending $\delta$ and such that
$\rho(u) = a$ and $\rho(v) = b$.
Since $\K \in \ECM$, there must exist $u', v' \in \K$ such that
$\tuple{u',v'} \in X$, $\delta u' = a$ and $\delta v' = b$.
\end{proof}


So, we want to show that any $\Kd \in \ECM$ defines the set $\Q$
of rational numbers: thus, no $\aleph_{1}$-saturated structure is in $\ECM$.\\

 Fix $\K^{\delta} \models \Td$.
Let $c \in K$ and define:
\begin{align*}
\Gamma_{c} &\coloneqq \set{\tuple{E(x) , E(cx)}: x \in K} \subseteq K^{2}\\
\delta_{c}&: \Gamma_{c} \to K^{2}\\
&\tuple{y,z} \mapsto \tuple{\delta y, \delta z}.
\end{align*}

By Lemma \ref{zariskidense} we have:

\begin{corollary}\label{lem:E-surjective}
Let  $\Kd \in \ECM,$  $\delta_{c}$ is surjective iff $\Gamma_{c}$ is Zariski dense.
\end{corollary}


\begin{lemma}\label{dense}
$\Gamma_{c}$ is not Zariski dense in $K^{2}$ iff $c \in \Q$.
\end{lemma}

It is then clear that the formula $\phi(t)$ saying  ``$\delta_{t}$ is not surjective''
defines $\Q$ in any $\Kd \in \ECM$.

\begin{proof}
If $c \in \Q$ then obviously $\Gamma_{c}$ is not Zariski dense.

\medskip For the converse, we assume that $c \in K \setminus \Q$, we need to show that $\Gamma_{c}$ is Zariski
dense.
If not, let $H$ be the Zariski closure of $\Gamma_{c}$ inside
$G \coloneqq ({\K^{*}})^{2}$.
Notice that $\Gamma_{c}$ is a subgroup of $G$, therefore $H$ is a Zariski closed
subgroup of $G$.
So by Corollary 3.2.15 in \cite{BG} there exist $m, n \in \Z$ not both zero, such that
$H \subseteq \set{\tuple{y,z} \in G: y^{n} \cdot z^{m} = 1}.$

Thus, for every $x \in K$,
\[
E(nx) \cdot E(mcx) = 1,
\]
and therefore $E(d x) = 1$, where $d \coloneqq n + cm$.
However, since $c \notin \Q$, $d \neq 0$, and therefore $E$
 is constant $1$, and this is a contradiction since the exponential function is not trivial.
\end{proof}

So, we have all ingredients to prove:

\begin{thm}\label{generale}
If $T$ expands a non-trivial exponential field, then $T$ does not admit generic derivations.
\end{thm}

\begin{proof}
It follows easily applying Corollary \ref{lem:E-surjective} and Lemma \ref{dense}.
\end{proof}

\begin{corollary}
$(\C, \exp),$ $(\mathbb B, E(x))$ where $\mathbb B$ is a Zilber field and $(\R,
\exp)$ do not admit generic derivations.
\end{corollary}



\section{Restricted exponential function}\label{altricasi}

Using similar idea we prove the case when there exists a definable restricted exponential function, that we denote by  $e: [-1, 1] \rightarrow \K^{*}$, where $\K$ is an extension of an ordered field, such that:

\[
e(0) = 1, \qquad
e(x - y) =\frac{ e(x)}{e(y)} \text { if }x, y, x-y \in [-1,1], \qquad
\exists z\, e(z) \neq 1.
\]

\begin{thm}
Assume that the function $e$ is not constant. 
Then, there exists an $L$-definable family $(\Gamma_{c}: c \in I)$ of subsets of $K^{2}$ with $I = [0, 1]$,
such that $\Gamma_{c}$ is not Zariski dense in $K^{2}$
iff $c \in \Q$.
\end{thm}

\begin{proof}
The set
\[
\Gamma_{c} \coloneqq \set{\tuple{e(x), e(cx)}: x \in [-1,1]}
\]
is uniformly $L$-definable.\\
We prove that $\Gamma_{c}$ is Zariski dense in $K^{2}$ iff $c \notin \Q$.\\
Assume that $c \in \Q$: then it is clear that $\Gamma_{c}$ is not Zariski dense.

\smallskip

Conversely, assume that $c \in I \setminus \Q$ but $\Gamma_{c}$ is not Zariski dense.
\Wlog, we may assume that $\K$ is $\omega$-saturated.
Let $S = \bigcup_{n \in \N} [-n, n]$; we extend $e(x)$ to all $S$ by
defining $e(n+x )=e^n \cdot e(x)$ where $n$ is an integer, $x \in [0, 1)$ 
and $e = e(1)$. It is easy to prove that $e : S \rightarrow \K^{*}$ is a group homomorphism. 
Let $\Gamma_{c, S} \coloneqq \set{\tuple{e(x), e(cx)}: x \in S}$:
notice that $\Gamma_{c,S}$ is a subgroup of $G \coloneqq (K^{*})^{2}$, and therefore
its Zariski closure $H$ is  a proper subgroup of~$G$.
Therefore, there exist (Corollary 3.2.15 in \cite{BG}) $m, n \in \Z$ not both zero such that
$\Gamma_{c,S} \subseteq \set{\tuple{y,z} \in K^{2}: y^{m} z^{n} = 1}$.
Thus, for every $x \in S$, $e(dx) = 1$, where $d \coloneqq m + nc$.
Notice that $d$ is finite and non-zero: thus, $e(x)$ is constant $1$, 
contradiction.
\end{proof}


Again by Lemma \ref{zariskidense} we have immediately the following:

\begin{corollary}
If $\delta_{c}$ is surjective then $c \notin \Q$.
If moreover $\Kd \in \ECM$ and $c \notin \Q$, then $\delta_{c}$ is surjective.
\end{corollary}

Using the previous results we can prove:

\begin{thm}
$(\R, \exp_{\upharpoonright [0, 1]})$ does not admit generic derivations.
\end{thm}

\section{Restricted sine}
\label{caso2}

Now we want to prove similar results when we consider the restricted $\sin x$. 
So, we consider
 $\sin(x), \cos(x) \in L$ unary function symbols and $\pi \in L$ a constant
symbol.  Assume that $T$ expands
the theory of ordered fields.
Let $\K \models T$, and 
denote by $D$ the interval $[0, 2\pi) \subset K$ with addition modulo $2\pi$; therefore,
$D$ is a definable Abelian group.
Denote by $\C$ the field  $\K + i \K$, and by $\C^{*}$ the multiplicative group of $\C$.

\begin{definition}
We say that $\K$ defines a restricted sine if $\pi > 0$, 
and the function
$F: D \to \C^{*}$, $x \mapsto \sin x + i \cos x$ is a group homomorphism which is not
constant.  
We say that $T$ defines a restricted sine if the above happens in all models
of~$T$.
\end{definition}

\begin{thm}\label{thm:sin}
If $T$ defines a restricted sine, 
then $T$ does not admit generic derivations.
\end{thm}

\begin{proof}
We can define
 $F(cx)$ uniformly in $x$ and $c$, as $c$ varies in $[0,1]$
and $x$ varies in $D$.
Thus, the set
\[
\Gamma_{c} \coloneqq \set{\tuple{F(x), F(cx)}: x \in D}
\]
is uniformly $L$-definable as $c \in [0,1]$.

Notice that $\Gamma_{c}$, as a subset of $K^4$, is never Zariski dense in~$K^4$. 
\begin{claim}\label{cl:C-dense}
$\Gamma_{c}$ is Zariski dense in $\C^{2}$ (with the Zariski topology of
$\C^{2}$, not the one of $\K^{4}$!) iff $c \notin \Q$.
\end{claim}
Assume that $c \in \Q$, then it is clear that $\Gamma_{c}$ is not Zariski dense.

Conversely, assume that $c \in [0,1] \setminus \Q$,
but, by contradiction, $\Gamma_{c}$ is not Zariski dense.
Notice that $\Gamma_{c}$ is a subgroup of $G \coloneqq ({\C^{*}})^{2}$;
therefore, its Zariski closure $H$ is a proper subgroup of~$G$.
Thus, reasoning as before, using  \cite[Corollary 3.2.15]{BG}, there exist
$m, n \in \Z$ not both zero such that, for every $x \in D$,
$F(dx) = 1$, where $d = m + nc$.
Notice that $d$ is  non-zero, and therefore $F$ is constant $1$, contradiction.

So $\Gamma_{c}$ is Zariski dense in $G$.

We claim that the real part of $\Gamma_{c},$ i.e. $\Gamma'_{c} := Re(\Gamma_{c})$ is Zariski dense in $K^{2}$ iff $c \notin \Q$.\\
It is easy to see that, if $c \in \Q$, $\Gamma'_{c}$  is Zariski dense in $K^2$, because
dimension of the projection of $\Gamma_{c}$ is preserved since the fibers are
finite. 
Now if we consider the derivation map $\delta_c : \Gamma'_c \rightarrow K^{2}$, we have that, if $c \in \Q$, since $\dim (\Gamma_c) \leq 1$, then $\delta_c$ is not surjective. If $c \not \in \Q$ then $\delta_c$ is surjective.
\end{proof}

\begin{corollary}
$(\R, \sin_{\upharpoonright [0, 2\pi]})$ does not admit generic derivations. 
\end{corollary}

\printbibliography

@book {BG,
    AUTHOR = {Bombieri, Enrico and Gubler, Walter},
     TITLE = {Heights in {D}iophantine geometry},
    SERIES = {New Mathematical Monographs},
    VOLUME = {4},
 PUBLISHER = {Cambridge University Press, Cambridge},
      YEAR = {2006},
     PAGES = {xvi+652},
      ISBN = {978-0-521-84615-8; 0-521-84615-3},
   MRCLASS = {11G50 (11-02 11G10 11G30 11J68 14G40)},
  MRNUMBER = {2216774},
MRREVIEWER = {Yuri\ Bilu},
       DOI = {10.1017/CBO9780511542879},
       URL = {https://doi.org/10.1017/CBO9780511542879},
}

@unpublished{point,
author={Point, Francoise and  Regnault, Nathalie}, 
title={Differential exponential topological fields},
eprinttype  = {arXiv},
eprint = {2007.14344v2},
year = {2023}
}

@article{FK,
eprinttype = {arXiv},
arxivId = {1905.07298},
eprint = {1905.07298},
author = {Fornasiero, Antongiulio and Kaplan, Elliot},
doi = {10.1142/S0219061321500070},
issn = {0219-0613},
journal = {Journal of Mathematical Logic},
month = {10},
title = {Generic derivations on o-minimal structures},
year = {2020}
}

@Unpublished{FT,
  author = 	 {Fornasiero, Antongiulio and Terzo, Giuseppina},
  title = 	 {Generic derivations on algebraically bounded structures},
  eprinttype = {arXiv},
  eprint = {2310.20511},
  year = 2024
}

@book{Sacks,
author = {Sacks, Gerald E},
title = {Saturated Model Theory},
publisher = {World Scientific},
year = {2009},
doi = {10.1142/6974},
edition   = {2nd},
OPTeprint = {https://www.worldscientific.com/doi/pdf/10.1142/6974}
}

@article {Robinson,
    AUTHOR = {Robinson, Abraham},
     TITLE = {On the concept of a differentially closed field},
   JOURNAL = {Bull. Res. Council Israel Sect. F},
  FJOURNAL = {Bulletin of the Research Council of Israel. Section F},
    VOLUME = {8F},
      YEAR = {1959},
     PAGES = {113--128},
   MRCLASS = {02.50 (12.80)},
  MRNUMBER = {125016},
MRREVIEWER = {R.\ C.\ Lyndon},
}

@book{ADH,
 title={Asymptotic Diﬀerential Algebra and Model Theory of Transseries}, ISBN={978-0-691-17543-0},
 doi = {10.1515/9781400885411},
 abstractNote={Asymptotic differential algebra seeks to understand the solutions of differential equations and their asymptotics from an algebraic point of view. The differential field of transseries plays a central role in the subject. Besides powers of the variable, these series may contain exponential and logarithmic terms. Over the last thirty years, transseries emerged variously as super-exact asymptotic expansions of return maps of analytic vector fields, in connection with Tarski's problem on the field of reals with exponentiation, and in mathematical physics. Their formal nature also makes them suitable for machine computations in computer algebra systems. This book validates the intuition that the differential field of transseries is a universal domain for asymptotic differential algebra. It does so by establishing in the realm of transseries a complete elimination theory for systems of algebraic differential equations with asymptotic side conditions. Beginning with background chapters on valuations and differential algebra, the book goes on to develop the basic theory of valued differential fields, including a notion of differential-henselianity. Next, H-fields are singled out among ordered valued differential fields to provide an algebraic setting for the common properties of Hardy fields and the differential field of transseries. The study of their extensions culminates in an analogue of the algebraic closure of a field: the Newton–Liouville closure of an H-field. This paves the way to a quantifier elimination with interesting consequences.},
 publisher={Princeton University Press},
 author={Aschenbrenner, Matthias and van den Dries, Lou and van der Hoeven, Joris},
 language={en} }

@article {Eklof,
    AUTHOR = {Eklof, Paul and Sabbagh, Gabriel},
     TITLE = {Model-completions and modules},
   JOURNAL = {Ann. Math. Logic},
  FJOURNAL = {Annals of Mathematical Logic},
    VOLUME = {2},
      YEAR = {1970/71},
    NUMBER = {3},
     PAGES = {251--295},
      ISSN = {0003-4843},
   MRCLASS = {02.50},
  MRNUMBER = {277372},
MRREVIEWER = {I.\ L.\ Gaal},
       DOI = {10.1016/0003-4843(71)90016-7},
       URL = {https://doi.org/10.1016/0003-4843(71)90016-7},
}

@article {zilber,
    AUTHOR = {Zilber, B.},
     TITLE = {Pseudo-exponentiation on algebraically closed fields of
              characteristic zero},
   JOURNAL = {Ann. Pure Appl. Logic},
  FJOURNAL = {Annals of Pure and Applied Logic},
    VOLUME = {132},
      YEAR = {2005},
    NUMBER = {1},
     PAGES = {67--95},
      ISSN = {0168-0072,1873-2461},
   MRCLASS = {03C60 (03C10 12L12)},
  MRNUMBER = {2102856},
MRREVIEWER = {Matthias\ Aschenbrenner},
       DOI = {10.1016/j.apal.2004.07.001},
       URL = {https://doi.org/10.1016/j.apal.2004.07.001},
}

@article {kirby,
    AUTHOR = {Haykazyan, Levon and Kirby, Jonathan},
     TITLE = {Existentially closed exponential fields},
   JOURNAL = {Israel J. Math.},
  FJOURNAL = {Israel Journal of Mathematics},
    VOLUME = {241},
      YEAR = {2021},
    NUMBER = {1},
     PAGES = {89--117},
      ISSN = {0021-2172,1565-8511},
   MRCLASS = {03C60 (03C45 12L12)},
  MRNUMBER = {4242146},
MRREVIEWER = {G.\ Cherlin},
       DOI = {10.1007/s11856-021-2089-1},
       URL = {https://doi.org/10.1007/s11856-021-2089-1},
}

\end{document}